\documentclass[12pt]{article}

\usepackage{graphicx,amsmath,amssymb,amsthm,amsfonts,latexsym,tikz,hyperref,subfig}
\usetikzlibrary{decorations.pathmorphing}
\usepackage[hmargin=.8in,vmargin=1in]{geometry}

 \DeclareMathOperator{\lrm}{lrm}

\newcommand{\fI}{{\mathfrak I}}

\newtheorem{thm}{Theorem}[section]
\newtheorem{prop}[thm]{Proposition}
\newtheorem{cor}[thm]{Corollary}
\newtheorem{lem}[thm]{Lemma}
\newtheorem{conj}[thm]{Conjecture}
\newtheorem{exa}[thm]{Example}

\newtheorem{defn}[thm]{Definition}

\newcommand{\ben}{\begin{enumerate}}
\newcommand{\een}{\end{enumerate}}
\newcommand{\ble}{\begin{lem}}
\newcommand{\ele}{\end{lem}}
\newcommand{\bth}{\begin{thm}}
\renewcommand{\eth}{\end{thm}}
\newcommand{\bpr}{\begin{prop}}
\newcommand{\epr}{\end{prop}}
\newcommand{\bco}{\begin{cor}}
\newcommand{\eco}{\end{cor}}
\newcommand{\bcon}{\begin{conj}}
\newcommand{\econ}{\end{conj}}
\newcommand{\bde}{\begin{defn}}
\newcommand{\ede}{\end{defn}}
\newcommand{\bex}{\begin{exa}}
\newcommand{\eex}{\end{exa}}
\newcommand{\barr}{\begin{array}}
\newcommand{\earr}{\end{array}}
\newcommand{\btab}{\begin{tabular}}
\newcommand{\etab}{\end{tabular}}
\newcommand{\beq}{\begin{equation}}
\newcommand{\eeq}{\end{equation}}
\newcommand{\bea}{\begin{eqnarray*}}
\newcommand{\eea}{\end{eqnarray*}}
\newcommand{\bal}{\begin{align*}}
\newcommand{\bce}{\begin{center}}
\newcommand{\ece}{\end{center}}
\newcommand{\bpi}{\begin{picture}}
\newcommand{\epi}{\end{picture}}
\newcommand{\bpp}{\begin{picture}}
\newcommand{\epp}{\end{picture}}
\newcommand{\bfi}{\begin{figure} \begin{center}}
\newcommand{\efi}{\end{center} \end{figure}}
\newcommand{\bprf}{\begin{proof}}
\newcommand{\eprf}{\end{proof}\medskip}

\newcommand{\bsl}{\begin{slide}{}}
\newcommand{\esl}{\end{slide}}
\newcommand{\bfr}{\begin{frame}}
\newcommand{\efr}{\end{frame}}

\newcommand{\pf}{{\bf Proof.\hspace{5pt}}}
\newcommand{\prf}{{\noindent\bf Proof.\hspace{5pt}}}
\newcommand{\hqed}{\hfill \qed}
\newcommand{\hqedm}{\hfill \qed \medskip}

\newcommand{\hso}[1]{\hspace{-1pt}}

\newcommand{\case}[4]{\left\{\barr{ll}#1&\mbox{#2}\\#3&\mbox{#4}\earr\right.}

\newcommand{\flf}[2]{\left\lfloor\frac{#1}{#2}\right\rfloor}
\newcommand{\cef}[2]{\left\lceil\frac{#1}{#2}\right\rceil}

\def\<{\langle}
\def\>{\rangle}

\newcommand{\ree}[1]{(\ref{#1})}

\newcommand{\si}{\sigma}

\newcommand{\1}{{\bf 1}}

\newcommand{\cA}{{\cal A}}

\newcommand{\fS}{{\mathfrak S}}

\DeclareMathOperator{\Av}{Av}

 \DeclareMathOperator{\des}{des}
\DeclareMathOperator{\Des}{Des}

 \DeclareMathOperator{\inv}{inv}

\DeclareMathOperator{\maj}{maj}

\begin{document}
\title{Inversion polynomials for $321$-avoiding permutations: addendum
}
\author{Szu-En Cheng\\[-5pt]
\small Department of Applied Mathematics, National University of Kaohsiung\\[-5pt]
\small Kaohsiung 811, Taiwan, ROC, {\tt chengszu@nuk.edu.tw}\\[5pt]
Sergi Elizalde\\[-5pt]
\small Department of Mathematics, Dartmouth College\\[-5pt]
\small Hanover, NH 03755-3551, USA, {\tt sergi.elizalde@dartmouth.edu}\\[5pt]
Anisse Kasraoui\\[-5pt]
\small Fakult\"at f\"ur Mathematik, Universit\"at Wien\\[-5pt]
\small Nordbergstra{\ss}e 15, A-1090 Vienna, Austria, {\tt anisse.kasraoui@univie.ac.at}\\[5pt]
Bruce E. Sagan\\[-5pt]
\small Department of Mathematics, Michigan State University,\\[-5pt]
\small East Lansing, MI 48824-1027, USA, {\tt sagan@math.msu.edu} }

\date{\today\\[10pt]
    \begin{flushleft}
                                           \\[5pt]
    \end{flushleft}}

\maketitle

\begin{abstract}
This addendum contains results about the inversion number and major index polynomials for permutations avoiding $321$ which did not fit well into the original paper.  In particular, we consider symmetry, unimodality, behavior modulo $2$, and signed enumeration.
\end{abstract}

\section{Basic definitions}

We recall the fundamental definitions from the original paper~\cite{ceks:ipt} for ease of reference.  

Call two sequences of distinct integers $\pi=a_1\ldots a_k$ and
$\si=b_1\ldots b_k$ {\em order isomorphic} whenever $a_i<a_j$ if and
only if $b_i<b_j$ for all $i,j$.  Let $\fS_n$ denote the symmetric
group of permutations of $[n]\stackrel{\rm def}{=}\{1,\ldots,n\}$.
Say that {\em $\si\in\fS_n$ contains $\pi\in\fS_k$ as a pattern} if
there is a subsequence $\si'$ of $\si$ order isomorphic to $\pi$.
If $\si$ contains no such subsequence then we say $\si$ {\em avoids}
$\pi$ and write $\Av_n(\pi)$ for the set of such $\si\in\fS_n$.

We will use a hash sign to denote cardinality.
Our generating functions will keep track of four statistics for $\si=b_1\ldots b_n\in\fS_n$:
\begin{itemize}
\item the {\em number of descents}
$$
\des\si=\#\Des\si
$$
where
$\Des\si=\{i\ |\ b_i>b_{i+1}\}$,
\item the {\em major index}
$$
\maj\si=\sum_{i\in\Des\si} i,
$$
\item the {\em inversion number}
$$
\inv\si=\#\{(i,j)\ |\ \text{$i<j$ and $b_i>b_j$}\},
$$
\item the number of {\em left-right maxima}
$$
\lrm\si=\#\{i\ |\ \text{$b_i>b_j$ for all $j<i$}\}.
$$
\end{itemize}
For $\Av(321)$ we will be studying the {\em major index polynomial}
$$
M_n(q,t)=\sum_{\si\in\Av_n(321)} q^{\maj\si} t^{\des\si}
$$
and the {\em inversion number polynomial}
$$
I_n(q,t)=\sum_{\si\in\Av_n(321)} q^{\inv\si} t^{\lrm\si}.
$$
Note that 
$$
M_n(1,1)=I_n(1,1)=\#\Av_n(321)=C_n
$$
where $C_n$ is the $n$th Catalan number.

\section{Symmetry, unimodality, and mod $2$ behavior of  $M_n(q,t)$}
\label{sum}

The coefficients of the polynomials $M_n(q,t)$ have various nice
properties which we now investigate. If $f(x)=\sum_k a_kx^k$ is a
polynomial in $x$ then we will use the notation
\bea
[x^k] f(x)&=&\mbox{coefficient of $x^k$ in $f(x)$}\\
&=&a_k.
\eea
Our main object of study in this section will be the polynomial
$$
A_{n,k}(q)=[t^k] M_n(q,t).
$$
In other words, $A_{n,k}(q)$ is the generating function for the
$\maj$ statistic over $\si\in\Av_n(321)$ having exactly $k$ descents. 

The first property which will concern us is symmetry. Consider a
polynomial
$$
f(x)=\sum_{i=r}^s a_i x^i
$$
where $a_r, a_s\neq0$.  Call $f(x)$ \emph{symmetric} if $a_i=a_j$
whenever $i+j=r+s$. \bth \label{sym} The polynomial $A_{n,k}(q)$ is
symmetric for all $n,k$. \eth \pf\ If $\si$ is counted by
$A_{n,k}(q)$ then $\des\si=k$.  Since $\si$ avoids $321$, it can not
have two consecutive descents and so the minimum value of $k$ is
$$
1+3+\cdots+(2k-1)=k^2
$$
and the maximum value is
$$
(n-1)+(n-3)+\cdots+(n-2k+1)=nk-k^2.
$$
So it suffices to show that for $0\le i\le nk$ we have
$a_i=a_{nk-i}$ where
$$
A_{n,k}(q)=\sum_i a_i q^i.
$$

Let $\cA_i$ be the permutations counted by $a_i$ and let $R_{\180}$
denote rotation of the diagram of $\si$ by $180$ degrees.  We claim
$R_{180}$ is a bijection between $\cA_i$ and $\cA_{nk-i}$ which will
complete the proof.  First of all, $R_{180}(321)=321$ and so
$\si$ avoids $321$ if and only if $R_{180}(\si)$ does so as well.  If $\si\in\cA_i$
then let $\Des\si=\{d_1,\ldots,d_k\}$ where $\sum_j d_j=i$.  It is
easy to see that $\Des R_{180}(\si)=\{n-d_1,\ldots,n-d_k\}$.  It
follows that $\maj R_{180}(\si)=nk-i$ and so $R_{180}(\si)\in
\cA_{nk-i}$.  Thus $R_{180}$ restricts to a well defined map from
$\cA_i$ to $\cA_{nk-i}$.  Since it is its own inverse, it is also a
bijection. \hqedm

Two other properties often studied for polynomials are unimodality
and log concavity.  The polynomial $f(x)=\sum_{i=0}^s a_i x^i$ is
{\em unimodal} if there is an index $r$ such that $a_0\le\ldots\le
a_r\ge\ldots\ge a_s$.  It is {\em log concave} if $a_i^2\ge a_{i-1}
a_{i+1}$ for all $0<i<s$.  If all the $a_i$ are positive, then log
concavity implies unimodality. \bcon The polynomial $A_{n,k}(q)$ is
unimodal for all $n,k$. \econ This conjecture has been checked by
computer for all $k< n\le 10$.  The corresponding conjecture for log
concavity is false, in particular, $A_{6,2}$ is not log concave.

The number-theoretic properties of the Catalan numbers have attracted
some interest.  Alter and Kubota~\cite{ak:ppp} determined the
highest power of a prime $p$ dividing $C_n$ using arithmetic means.
Deutsch and Sagan~\cite{ds:ccm}  gave a proof of this result using
group actions for the special case $p=2$.  Just considering parity,
one gets the nice result that $C_n$ is odd if and only if $n=2^m-1$
for some nonnegative integer $m$.  Dokos et al. proved the following
refinement of the ``if" direction of this statement. \bth[Dokos et
al.~\cite{ddjss:pps}] Suppose $n=2^m-1$ for some $m\ge0$.  Then
$$
\hspace*{130pt}[q^k]I_n(q,1)=\case{1}{if $k=0$}{\text{an even
integer}}{if $k\ge 1$.\hspace*{130pt}\qed}
$$
\eth In the same paper, the following statement was made as a
conjecture which has now been proved by Killpatrick.

\bth[Killpatrick~\cite{kil:pcc}] Suppose $n=2^m-1$ for some $m\ge0$.
Then
$$
\hspace*{130pt}[q^k]M_n(q,1)=\case{1}{if $k=0$}{\text{an even
integer}}{if $k\ge 1$.\hspace*{130pt}\qed}
$$
\eth

\bfi
\begin{tikzpicture}
\draw (0,0) -- (2,0) -- (2,2) -- (0,2) -- (0,0); \draw (2/3,0) --
(2/3,2/3) -- (0,2/3); \draw (2,2/3) -- (4/3, 2/3) -- (4/3,2);; \draw
(2,4/3) -- (2/3,4/3) -- (2/3,2); \node at (1/3,1/3){$\si_1$}; \node
at (1,5/3){$\si_2$}; \node at (5/3,1){$\si_3$};
\end{tikzpicture}
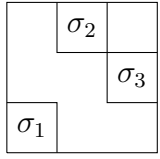
\caption{The diagram of $132[\si_1,\si_2,\si_3]$} \label{132} \efi

We wish to prove a third theorem of this type.  To do so, we will
need the notion inflation for permutations.  Given a permutation
$\pi=a_1\ldots a_n\in\fS_n$ and permutations $\si_1,\ldots,\si_n$,
the {\em inflation} of $\pi$ by the $\si_i$, written
$\pi[\si_1,\ldots,\si_n]$, is the permutation whose diagram is
obtained from the diagram of $\pi$ by replacing the dot $(i,a_i)$ by
a copy of $\si_i$ for $1\le i\le n$.  By way of example,
Figure~\ref{132} shows a schematic diagram of an inflation of the
form $132[\si_1,\si_2,\si_3]$.  More specifically,
$132[21,1,312]=216534$.

\bth Suppose $n=2^m-1$ for some  $m\ge0$.  Then
$$
[t^k]M_n(1,t)=A_{n,k}(1)=\case{1}{if $k=0$,}{\mbox{an even
integer}}{if $k\ge1$.}
$$
\eth \pf\ We have  $A_{n,0}(q)=1$ since $\si=12\ldots n$ is the only
permutation without descents and it avoids $321$.

If $k\ge1$, and $A_{n,k}(q)$ has an even number of terms then
$A_{n,k}(1)$ must be even because it is a symmetric polynomial by
Theorem~\ref{sym}. By the same token, if $A_{n,k}(q)$ has an odd
number of terms, then $A_{n,k}(1)$ has the same parity as its middle
term.  Now  consider $R_{180}$ acting on the elements of
$\cA_{nk/2}$ as in the previous proof.  Note that since $n$ is odd,
$k$ must be even.  Furthermore, this action partitions $\cA_{nk/2}$
into orbits of size one and two.  So it suffices to show that there
are an even number of fixed points.  If $\si$ is fixed then its
diagram must contain the center, $c$,  of the square since this is a
fixed point of $R_{180}$.  Also, the NW and SE quadrants of $\si$
with respect to $c$ must be empty, since otherwise they both must
contain dots (as one is taken to the other by $R_{180}$) and
together with $c$ this forms a $321$.   For the same reason, the SW
quadrant of $\si$ determines the NE one. Thus, the fixed points are
exactly the inflations of the form $\si=123[\tau,1,R_{180}(\tau)]$
where $\tau\in\Av_{2^{m-1}-1}(321)$ has $k/2$ descents.  By
induction on $m$ we have that  the number of such $\tau$, and hence the number of
such $\si$, is even. \hqedm

%
%

\section{Refined sign-enumeration of $321$-avoiding permutations}
\label{rse}

Simion and Schmidt~\cite{ss:rp} considered the signed
enumeration of various permutation classes of the form
$\sum_{\si\in\Av_n(\pi)}(-1)^{\inv\si}$.  In this section we will
rederive their theorem for $\Av_n(321)$ using a result from~\cite{ceks:ipt} about continued fractions.  In addition, we will  provide a more refined signed
enumeration which also keeps track of the $\lrm$ statistic.  We
should note that Reifegerste~\cite{rei:rsb} also has a refinement
which takes into account the length of the longest increasing
subsequence of $\si$.

We will use the following notation for continued
fractions \beq \label{cf} F=
\frac{a_1|}{|b_1}\pm\frac{a_2|}{|b_2}\pm\frac{a_3|}{|b_3}\pm\cdots
\quad = \quad \cfrac{a_1}{b_1\pm \cfrac{a_2}{b_2\pm
\cfrac{a_3}{b_3\pm\dotsb }}}. \eeq 
Now consider the generating
function $C(z)=\sum_{n\ge0} C_n z^n$. It is well known that $C(z)$
satisfies the functional equation $C(z)=1+zC(z)^2$. Rewriting this
as $C(z)=1/(1-z C(z))$ and iteratively substituting for $C(z)$, we
obtain the also well-known continued fraction
\beq
\label{C(z)}
C(z)=\frac{1|}{|1}-\frac{z|}{|1}-\frac{z|}{|1}-\frac{z|}{|1}-\frac{z|}{|1}-\frac{z|}{|1}-\cdots.
\eeq

To return to our context, consider
$$
{\mathfrak I}(q,t;z)=\sum_{n\ge0} I_n(q,t) z^n.
$$
We have a continued fraction expansion for this power series generalizing the one in equation~\ree{C(z)}.
\bco[Cheng et al.~\cite{ceks:ipt}] The generating function $\fI(q,t;z)$
has  continued fraction expansion 
\beq \label{sti}
\fI(q,t;z) =\frac{1|}{|1}-\frac{tz|}{|1}-\frac{qz|}{|1}
 -\frac{t qz|}{|1}-\frac{q^2 z|}{|1}-\frac{tq^2 z|}{|1}-\frac{q^3 z|}{|1}-\frac{tq^3 z|}{|1}
-\frac{q^4 z|}{|1}-\frac{tq^4 z|}{|1}-\cdots \eeq \eco

We will also need the following well-known result.
\bth[Jones and Thron~\cite{jt:cf}] \label{jt:thm} We have \bea
\frac{a_1|}{|1}+\frac{a_2|}{|1}+\frac{a_3|}{|1}+\cdots &=&
\frac{a_1|}{|1+a_2}-\frac{a_2a_3|}{|1+a_3+a_4}-\frac{a_4a_5|}{|1+a_5+a_6}-\frac{a_6a_7|}{|1+a_7+a_8}+\cdots\\[10pt]
&=&
a_1-\frac{a_1a_2|}{|1+a_2+a_3}-\frac{a_3a_4|}{|1+a_4+a_5}-\frac{a_5a_6|}{|1+a_6+a_7}
-\frac{a_7a_8|}{|1+a_8+a_9}-\cdots \eea 
where the second and third continued fractions are called the even and odd parts,
respectively,  of the first continued fraction.\hqed \eth

Now plug $q=-1$ and $t=1$ into the continued fraction~\ree{sti} to
obtain
$$
\fI(-1,1;z)=\frac{1|}{|1}-\frac{z|}{|1}+\frac{z|}{|1}
 +\frac{z|}{|1}-\frac{z|}{|1}-\frac{z|}{|1}+\frac{z|}{|1}
 +\frac{z|}{|1}-\frac{z|}{|1}-\frac{z|}{|1}+\frac{z|}{|1}
 +\frac{z|}{|1}\cdots.
$$
Using Theorem~\ref{jt:thm} to extract the odd part of this expansion
gives
$$
\fI(-1,1;z)=1+\frac{z|}{|1}-\frac{z^2|}{|1}-\frac{z^2|}{|1}-\frac{z^2|}{|1}-\frac{z^2|}{|1}-\frac{z^2|}{|1}-\cdots.
$$
Comparing this to the continued fraction for $C(z)$ in~\ree{C(z)}, we see that
$$
\fI(-1,1;z)=1+zC(z^2).
$$
Taking the coefficient of $z^n$ on both sides yields the following
result. \bth[Simion and Schmidt~\cite{ss:rp}] For all $n\ge1$, we
have
$$
 I_{2n}(-1,1)=\sum_{\sigma\in\Av_{2n}(321)}(-1)^{\inv\si}=0\quad\text{and}\quad
 I_{2n+1}(-1,1)=\sum_{\sigma\in\Av_{2n+1}(321)}(-1)^{\inv\si}=C_{n}.
$$  \hqed \eth

Since our refined sign-enumeration will involve the parameter
$\lrm$, we recall (but will not use) the folklore result that the enumerating
polynomial of $\Av_n(321)$ according to the $\lrm$ statistic is the
$n$th {\em Narayana polynomial}, i.e.,
 $$
 I_{n}(1,t)=\sum_{\sigma\in\Av_{n}(321)}t^{\lrm\si}=\sum_{k=1}^n N_{n,k} t^k,
 $$
where the {\em Narayana number} $N_{n,k}$ is given by
$N_{n,k}=\frac{1}{n}{n\choose k}{n\choose k-1}$ for $n\ge k\ge1$.

\bth \label{thm:sign-enumeration} For all $n\geq 1$,
 \beq \label{eq:sign-enum}
  I_{n}(-1,t)=\sum_{\sigma\in\Av_{n}(321)}(-1)^{\inv\si}t^{\lrm\si}=\sum_{k=1}^{n}(-1)^{n-k} s_{n,k}
  t^{k}
 \eeq
where $s_{n,k}$ is defined for $n\ge k\ge1$ by
$$s_{n,k}={\flf{n-1}{2} \choose \flf{k-1}{2}} {\cef{n-1}{2} \choose\cef{k-1}{2}}.$$
 Moreover,
 \begin{eqnarray}
 I_{2n}(-1,t)&=&(t-1)I_{2n-1}(-1,t),\label{eq:rec1-sign}\\
 (n+1)I_{2n+1}(-1,t)&=&2((1+t^2)n-t)I_{2n-1}(-1,t)-(1-t^2)^2(n-1)I_{2n-3}(-1,t).\label{eq:rec2-sign}
 \end{eqnarray}
\eth

\prf Let $\fI(t;z)$ and $\fI_{odd}(t;z)$ be the power series defined
as
$$
\fI(t;z)=\sum_{n\ge0} I_{n}(-1,t)z^n\quad\text{and}\quad
\fI_{odd}(t;z)=\sum_{n\ge0} I_{2n+1}(-1,t)z^{n}.
$$
By equation~\ree{sti}, we have
\bea \fI(t;z) &=&
\frac{1|}{|1}-\frac{tz|}{|1}+\frac{z|}{|1}
 +\frac{tz|}{|1}-\frac{z|}{|1}-\frac{tz|}{|1}+\frac{z|}{|1}+\frac{tz|}{|1}
-\frac{z|}{|1}-\frac{tz|}{|1}+\frac{z|}{|1}+\frac{tz|}{|1}-\cdots
\\[10pt]
 &=&
\cfrac{1}{1- \cfrac{tz}{1+ \cfrac{z}{1+ \cfrac{tz}{1-z\fI(t;z)}}}}
\eea which, after simplification, leads to the functional equation
 $$
 (1+z-tz)z\,\fI(t;z)^2 -(1+2z+z^2-t^2z^2)\,\fI(t;z)+(1+z+tz)=0.
$$ Solving this quadratic equation, we obtain
 \beq\label{eq:GF-sign}
 \fI(t;z)=\frac{1+2z+(1-t^2)z^2-\sqrt{
 1-2(1+t^2)z^2+(1-t^2)^2z^4}}{2z(1+z-tz)}.
 \eeq
Noticing that
$\fI_{odd}(t;z)=(\fI(t;\sqrt{z})-\fI(t;-\sqrt{z}))/2\sqrt{z}$ and
using~\eqref{eq:GF-sign}, we obtain after a routine computation
 \beq\label{eq:GF-sign-odd}
 \fI_{odd}(t;z)=\frac{1-(1-t)^2z-\sqrt{1-2(1+t^2)z+(1-t^2)^2z^2}}{2z(1-(1-t)^2z)}.
 \eeq

 It follows from~\eqref{eq:GF-sign} that
$$
 (1+z-tz)\fI(t;z)+(1-z+tz)\fI(t;-z)=2.
$$
Extracting the coefficient of $z^{2n}$ on both sides of the last
equality, we obtain~\eqref{eq:rec1-sign}.
Using~\eqref{eq:GF-sign-odd}, it is easily checked  that $\fI_{odd}$
satisfies the differential equation
$$
 z\big(1-2(1+t^2)z+(1-t^2)^2z^2\big)\,\fI_{odd}'(t;z)+\big(1-2(1-t+t^2)z+(1-t^2)^2z^2)\big)\,\fI_{odd}(t;z)-t=0,
$$
 where $\fI_{odd}'(t;z)$ is the derivative with respect to $z$.
 Extracting the coefficient of $z^{n}$ on both sides of the last
equality, we obtain~\eqref{eq:rec2-sign}.\\

We now turn our attention to~\eqref{eq:sign-enum}. Clearly, we have
 \beq\label{eq:coeff-signenum-a}
 [t^{2k+1}]I_{2n+1}(-1,t)=[t^k z^n]\frac{\fI_{odd}(\sqrt{t};z)-\fI_{odd}(-\sqrt{t};z)}{2\sqrt{t}}
                         =[t^{k}][z^n] \frac{1}{\sqrt{1-2(1+t)z+(1-t)^2z^2}}
 \eeq
where the last equality follows from~\eqref{eq:GF-sign-odd}. Using
the Lagrange inversion formula, one can show that
 \beq
 [z^n]\frac{1}{\sqrt{1-2(1+t)z+(1-t)^2z^2}}=[x^n](1+(1+t)x+tx^2)^n.
 \eeq
Combining~\eqref{eq:coeff-signenum-a} with the above relation, we
obtain
 \beq\label{eq:coeff-signenum-for-a}
 [t^{2k+1}]I_{2n+1}(-1,t)=[x^n][t^{k}](1+(1+t)x+tx^2)^n=[x^n]\binom{n}{k}x^k(1+x)^n=\binom{n}{k}^2.
 \eeq
Similarly, we have
 \bea
 [t^{2k}]I_{2n+1}(-1,t)
  =[t^{k}z^n]\frac{\fI_{odd}(\sqrt{t};z)+\fI_{odd}(-\sqrt{t};z)}{2}
  =[t^{k}] [z^n]\frac{1}{2z}-\frac{1-z-tz}{2z\sqrt{1-2(1+t)z+(1-t)^2z^2}}.
 \eea
This, combined first with~\eqref{eq:coeff-signenum-a} and
then~\eqref{eq:coeff-signenum-for-a}, yields
$$
 [t^{2k}]I_{2n+1}(-1,t)=\frac{1}{2}\big([t^{2k+1}]I_{2n+1}(t,-1)+[t^{2k-1}]I_{2n+1}(t,-1)-[t^{2k+1}]I_{2n+3}(t,-1)\big)
  =-\binom{n}{k-1}\binom{n}{k}.
$$
This proves that~\eqref{eq:sign-enum} is true when $n$ is odd.
Combining this with~\ree{eq:rec1-sign} shows that the formula also holds when $n$ is even.
 \hqedm

To see why the previous result implies the one of Simion and Schmidt, plug $t=1$ into the equations for $I_{2n}(-1,t)$ and $I_{2n+1}(-1,t)$.  In the former case we immediately get $I_{2n}(-1,1)=0$ because of the factor of $t-1$ on the right.  In the latter, we get the equation $(n+1) I_{2n+1}(-1,1)=2(2n-1) I_{2n-1}(-1,1)$.  The fact that $I_{2n+1}(-1,1)=C_n$ now follows easily by induction.

Finally, it is interesting to note that the
numbers $s_{n,k}$ which arise in the signed enumeration of
$\Av_n(321)$ have a nice combinatorial interpretation.  Recall that
{\em symmetric} Dyck paths are those $P=s_1\ldots s_{2n}$ which are
the same read forwards as read backwards. The following result
appears in Sloane's Encyclopedia~\cite{slo:oei}: For $n\geq
k\geq1$, the number $s_{n,k}$ is equal to the number of symmetric
Dyck paths of semilength $n$ with $k$ peaks.


\begin{thebibliography}{1}

\bibitem{ak:ppp}
Ronald Alter and K.~K. Kubota.
\newblock Prime and prime power divisibility of {C}atalan numbers.
\newblock {\em J. Combinatorial Theory Ser. A}, 15:243--256, 1973.

\bibitem{ceks:ipt}
Szu-En Cheng, Sergi Elizalde, Anisse Kasraoui, and Bruce~E. Sagan.
\newblock Inversion polynomials for $321$-avoiding permutations.
\newblock Preprint {\texttt{arXiv:1112.6014}}.

\bibitem{ds:ccm}
Emeric Deutsch and Bruce~E. Sagan.
\newblock Congruences for {C}atalan and {M}otzkin numbers and related
  sequences.
\newblock {\em J. Number Theory}, 117(1):191--215, 2006.

\bibitem{ddjss:pps}
Theodore Dokos, Tim Dwyer, Bryan~P. Johnson, Bruce~E. Sagan, and Kimberly
  Selsor.
\newblock Permutation patterns and statistics.
\newblock {\em Discrete Math.}, 312(18):2760--2775, 2012.

\bibitem{jt:cf}
William~B. Jones and Wolfgang~J. Thron.
\newblock {\em Continued fractions}, volume~11 of {\em Encyclopedia of
  Mathematics and its Applications}.
\newblock Addison-Wesley Publishing Co., Reading, Mass., 1980.
\newblock Analytic theory and applications, With a foreword by Felix E.
  Browder, With an introduction by Peter Henrici.

\bibitem{kil:pcc}
Kendra Killpatrick.
\newblock On the parity of certain coefficients for a {$q$}-analogue of the
  {C}atalan numbers.
\newblock {\em Electron. J. Combin.}, 19(4):Paper 27, 7, 2012.

\bibitem{rei:rsb}
Astrid Reifegerste.
\newblock Refined sign-balance on 321-avoiding permutations.
\newblock {\em European J. Combin.}, 26(6):1009--1018, 2005.

\bibitem{ss:rp}
Rodica Simion and Frank~W. Schmidt.
\newblock Restricted permutations.
\newblock {\em European J. Combin.}, 6(4):383--406, 1985.

\bibitem{slo:oei}
N.~J.~A. Sloane.
\newblock The on-line encyclopedia of integer sequences.
\newblock {\em Notices Amer. Math. Soc.}, 50(8):912--915, 2003.

\end{thebibliography}
\end{document}